\documentclass[12pt]{article}

\begin{document}

\title{Equivalence of the Bergman and Teichm\"uller metrics on Teichm\"uller spaces}
\author{Bo-Yong Chen}
\date{}
\maketitle

\section{Introduction}
Let $R$ be a compact Riemann surface with genus $g>1$. Denote by
Teich$(R)$ the Teichm\"uller space of $R$. There are two canonical
invariant metrics on Teich$(R)$, namely, the Teichm\"uller and
Weil-Petersson metrics. By Bers embedding one can regard
Teich$(R)$ as a bounded domain of holomorphy in ${\bf C}^{3g-3}$.
Hence it carries four classical invariant metrics: the
Carath$\acute{\rm e}$odory, Bergman, Kobayashi and
K\"ahler-Einstein metrics. Royden's theorem [12] states that the
Teichm\"uller metric coincides with the Kobayashi metric. The
Weil-Petersson metric is incomplete (cf. [6], [15]), while all the
other metrics are complete since the Carath$\acute{\rm e}$odory
metric is (cf. [3]).  Recently, McMullen [7] introduced a new
invariant metric $g_{1/l}$ equivalent to the Teichm\"uller metric
in order to prove that the moduli space of Riemann surfaces is
K\"ahler hyperbolic. By using the McMullen metric, Liu, Sun and
Yau [5] proved the equivalence of the Teichm\"uller and
K\"ahler-Einstein metrics. In this paper, we use McMullen's
$g_{1/l}$ metric and Takhtajan-Teo's K\"ahler potential of the
Weil-Petersson metric [13] together with the classical
$L^2-$estimate to show the following

\bigskip

\textbf{Theorem 1.1.} \emph{The Bergman and Teichm\"uller metrics
are equivalent on ${\rm Teich}(R)$.}

\bigskip

Throughout the paper, $A=O(B)$ means $A\le CB$ and equivalence
$A\asymp B$ means $\frac1C B\le A\le CB$ where $C>0$ is a uniform
constant on Teich$(R)$.

It is also interesting to invest the boundary behavior of the
Bergman kernel and metric if we regard Teich$(R)$ as a bounded
domain in ${\bf C}^{3g-3}$. Let $K_T$ denote the Bergman kernel
function, let $dist_B$ be the Bergman distance and $\delta_T$ be
the Euclidean boundary distance. We have

\bigskip

\textbf{Theorem 1.2.} \emph{$(i)\ K_T\ge
C(\delta_T|\log\delta_T|)^{-2}$; $(ii)$ Given $X_0\in {\rm
Teich}(R)$, we have $dist_B(X_0,\cdot)\ge C|\log\delta_T|$.}

\bigskip

\textbf{Remark.} Ohsawa [10] has showed $K_T(X)\rightarrow \infty$
as $X\rightarrow \partial {\rm Teich}(R)$.

\section{A review of Teichm\"uller theory}
In this section, we review some basic definitions in the
Teichm\"uller theory, for more detail, see [7].

\bigskip

A Riemann surface $R$ is called \emph{hyperbolic} if its universal
covering is the upper half plane $H$. The Poincar$\acute{{\rm e
}}$ metric $|dz|/{\rm Im\,}z$ on $H$ descends to a complete metric
on $R$ with constant curvature $-1$, which is called the
\emph{hyperbolic metric}.

\bigskip

Let $R$ be a hyperbolic Riemann surface. A Riemann surface $X$ is
\emph{marked} by $R$ if it is equipped with a qusiconformal
homeomorphism $f: R \rightarrow X$. Two marked surfaces
$(X_1,f_1),\,(X_2,f_2)$ are \emph{equivalent} if $f_2\circ
f_1^{-1}$ is homotopic to a conformal mapping of $X_1$ onto $X_2$.
We call the set of all such equivalence classes $[X,f]$ the
\emph{Teichm\"uller space} of $R$ and denote it by Teich$(R)$. The
\emph{Teichm\"uller distance} between two points
$p_i=[X_i,f_i],\,i=1,2$ in Teich$(R)$ is defined by
$$
d_T(p_1,p_2)=\frac12 \inf\log K(h)
$$
where $h$ is taken over all quasiconformal mappings of $X_1$ onto
$X_2$ which are homotopic to $f_2\circ f_1^{-1}$ and $K(h)\ge 1$
denotes the maximal dilatation of $h$. The Teichm\"uller space is
topologically a cell.

Given $X\in {\rm Teich}(R)$, let $Q(X)$ denote the Banach space of
holomorphic quadratic differentials $\phi=\phi(z)dz^2$ on $X$ with
$$
\|\phi\|_T=\int_X |\phi|<\infty.
$$
Let $B(X)$ be the space of $L^\infty$ measurable Beltrami
differentials $\mu=\mu(z)d\bar{z}/dz$ on $X$. A tangent vector
$v\in T_X{\rm Teich}(R)$ is represented by a $\mu\in B(X)$ and its
\emph{Teichm\"uller norm} is given by
$$
\|\mu\|_T=\sup\left\{{\rm Re\,}\int_X
\phi(z)\mu(z)dzd\bar{z}:\phi\in Q(X),\ \|\phi\|_T=1 \right\}.
$$
We have the isomorphism
$$
T_X{\rm Teich}(R)\cong B(X)/Q(X)^{\perp}
$$
and $\|\mu\|_T$ gives infinitesimal form of the Teichm\"uller
distance.

\bigskip

A \emph{projective structure} on $X$ is a subatlas charts with
M\"obius transformations as transition functions. The space of all
projective surfaces marked by $R$ is a complex manifold fibering
over Teich$(R)$, which will be denoted by Proj$(R)$. By Fuchsian
Uniformization, there is a canonical section $\sigma_F:{\rm
Teich}(R)\rightarrow {\rm Proj}(R)$. Each fiber ${\rm Proj}_X(R)$
over $X\in {\rm Teich}(R)$ is an affine space modeled on the
Banach space $P(X)$ of all holomorphic quadratic differentials on
$X$ with
$$
\sup_X\rho^{-2}|\phi|<\infty.
$$
Teich$(R)$ has a complexification defined by
$$
QF(R)={\rm Teich}(R)\times {\rm Teich}(\overline{R})
$$
where $\overline{R}$ is the complex conjugate of $R$. The
real-analytic map $\sigma_F$ naturally induces a holomorphic map
$$
\sigma:QF(R)\rightarrow {\rm Proj}(R)\times {\rm
Proj}(\overline{R}).
$$
Denote by
$\sigma(X,Y)=(\sigma_{QF}(X,Y),\overline{\sigma}_{QF}(X,Y))$. The
\emph{Bers embedding} $\beta_Y:{\rm Teich}(R)\rightarrow P(Y)$ is
given by
$$
\beta_Y(X)=\overline{\sigma}_{QF}(X,Y)-\sigma_F(Y).
$$
One has the following well-known theorem

\bigskip

 \textbf{Theorem 2.1.} \emph{The Bers embedding maps ${\rm Teich}(R)$ to a bounded domain in $P(Y)$ which
 is contained in the ball with radius $3/2$.}

\section{Weil-Petersson metric}

Let $R$ be a compact Riemann surface of genus $g>1$. The
\emph{Weil-Petersson norm} on the cotangent space $Q(X)\cong
T_X^\ast {\rm Teich}(R)$ is defined by
$$
\|\phi\|^2_{WP}=\int_X \rho(z)^{-2}|\phi(z)|^2 |dz|^2.
$$
By duality, one gets a Riemann metric $g_{WP}$ on the tangent
space of Teich$(R)$. Furthermore, it is a non-complete K\"ahler
metric of negative sectional curvature (cf. [6], [14]--[16]). It
follows from the Cauchy-Schwarz inequality that
$$
\|v\|_{WP}\le 2\sqrt{\pi(g-1)}\|v\|_T
$$
holds for any tangent vector $v$ on Teich$(R)$. Recall that
$$
\beta_X:{\rm Teich}(\overline{R})\rightarrow Q(X)\cong T_X^\ast
{\rm Teich}(R).
$$
It was shown by McMullen that for any fixed $Y\in {\rm
Teich}(\overline{R})$, the $1-$form $\beta_X(Y)$ is bounded in the
Teichm\"uller and Weil-Petersson metrics and satisfies
$d\beta_X(Y)=i\omega_{WP}$ where $\omega_{WP}$ is the K\"ahler
form of $g_{WP}$ (cf. Theorem 1.5 in [7]). Moreover, Takhtajan and
Teo found a real-analytic function $S_Y$ on Teich$(R)$ coming from
the \emph{Liouville action} in string theory such that
$$
-\beta_X(Y)=\sigma_F(X)-\sigma_{QF}(X,Y)=\frac12 \partial S_Y
$$
(cf. Corollary 4.1 in [13]), which implies $-\frac12 S_Y$ is a
K\"ahler potential for the Weil-Petersson metric with
\begin{equation}
\partial\bar{\partial} (-S_Y)\ge C\partial S_Y\bar{\partial} S_Y
\end{equation}
for suitable constant $C>0$.

\bigskip

Given a hyperbolic geodesic $\gamma$ on $R$, let $l_\gamma(X)$
denote the hyperbolic length of the corresponding geodesic on
$X\in {\rm Teich}(R)$. The length function is very useful in the
Teichm\"uller theory. For instance, it relates the Teichm\"uller
metric as follows
\begin{equation}
\|\partial \log l_\gamma \|_T\le 2
\end{equation}
(cf. Theorem 4.2 in [7]). Let ${\rm Log}:{\bf R}_+\rightarrow
[0,\infty)$ be a smooth function such that
$$
{\rm Log}(t)=\left\{
\begin{array}{ll}
\log t & {\rm if\ }t\ge 2\\
0 & {\rm if\ }t\le 1.
\end{array}
\right.
$$
McMullen [7] defined a new invariant K\"ahler metric by
$$
g_{1/l}=g_{WP}-\delta \sum_{l_\gamma(X)<\epsilon}
\partial\bar{\partial} {\rm Log}\frac{\epsilon}{l_\gamma}
$$
where the sum is over primitive short geodesics $\gamma$ on $X$;
at most $3g-3$ terms occur in the sum.

\bigskip

\textbf{Theorem.} (cf. Theorem 5.1 in [7]) \emph{For all
$\epsilon>0$ sufficiently small, there exists a $\delta>0$ such
that $g_{1/l}$ is equivalent to the Teichm\"uller metric.}

\bigskip

Set
$$
\psi=-\frac{S_Y}2-\delta \sum_{l_\gamma(X)<\epsilon}  {\rm
Log}\frac{\epsilon}{l_\gamma}.
$$

\bigskip

\textbf{Proposition 3.1.} \emph{There is a constant $C>0$ such
that}
\begin{equation}
g_{1/l}=\partial\bar{\partial}\psi\ge C\partial
\psi\bar{\partial}\psi.
\end{equation}

\bigskip

\emph{Proof.} It suffices to show the inequality in (3). For any
$\epsilon/2<l_\gamma(X)<\epsilon$, one has
$$
\|\partial {\rm Log}( \epsilon/l_\gamma)\|_T \le \sup_{t\in
[\epsilon/2,\epsilon]} |{\rm Log}'(t)|\cdot
\frac{\epsilon}{l_\gamma} \|\partial \log l_\gamma\|_T=O(1)
$$
by (2). By (1), (2) and the above theorem, the desired inequality
follows immediately from the Cauchy-Schwarz inequality.

\section{Proofs of Theorems 1.1 and 1.2}
Let $M$ be a complex manifold of dimension $m$. Let ${\cal H}_1$
denote the space of holomorphic $m-$forms $s$ on $M$ such that
$$
\|s\|^2_2=\left |\int_M s\wedge \bar{s} \right|\le 1.
$$
The \emph{Bergman kernel} on $M$ is defined by
\begin{equation}
K_M(z)=\sup\{s\wedge \bar{s}(z):s\in {\cal H}_1 \}
\end{equation}
where $s_1\wedge \bar{s}_1(z)\le  s_2\wedge \bar{s}_2(z)$ means
the ratio of the left and right sides is bounded by 1. If $K_M$ is
nowhere vanishing on $M$, one can define the \emph{Bergman metric}
by $ds^2_M:=\partial\bar{\partial}\log K^\ast_M$ where
$$
K_M=K_M^\ast dz_1\wedge \cdots\wedge dz_m\wedge d\bar{z}_1\wedge
\cdots \wedge d\bar{z}_m
$$
in local coordinates (note that the definition of $ds_M^2$ does
not depend on the choice of coordinates hence is globally
defined). It has the following extreme property:
\begin{eqnarray}
ds^2_M(z;v)  =  \frac1{K_M^\ast(z)}\sup\{ && |\partial
s^\ast(v)|^2(z): s\in {\cal H }_1,\ s(z)=0,\nonumber\\
&& s=s^\ast dz_1\wedge \cdots \wedge dz_n\ \ \}
\end{eqnarray}
for all $v\in T_z M$.

\bigskip

By Royden's theorem [12], given $X_0\in {\rm Teich}(R)$, there is
an embedded polydisk
$$
\iota:(\Delta^{3g-3},0)\rightarrow ({\rm Teich}(R),X_0)
$$
such that the Teichm\"uller(=Kobayashi) and Euclidean metrics are
equivalent on $\Delta^{3g-3}$. For any $s\in {\cal H}_1$ with
$s=s^\ast dz_1\wedge\cdots \wedge dz_{3g-3}$ in $\Delta^{3g-3}$,
we obtain from Cauchy's estimate that
\begin{equation}
|\partial^\alpha s^\ast/\partial z^\alpha (0)| = O(1)
\end{equation}
holds for any multi-indices $\alpha$. Let $\psi$ be as in section
3. By McMullen's theorem and Proposition 3.1, one has
\begin{eqnarray*}
|\psi(X)-\psi(X_0)| & \le & \left(\sup_{{\rm Teich}(R)}\|d
\psi\|_T\right)d_T(X_0,X)\\
& \le & O(d_T(X_0,X)).
\end{eqnarray*}
Thus there exists a constant $c_0>0$ independent of $X_0$ such
that
\begin{equation}
\iota (\Delta^{3g-3})\subset \{X\in {\rm
Teich}(R):|\psi(X)-\psi(X_0)|<c_0\}.
\end{equation}
Set
$$
\lambda(X)=-e^{-\frac{C}{2}(\psi(X)-\psi(X_0))}.
$$
By (3), one has
\begin{eqnarray*}
\partial\bar{\partial} \lambda & = & -\frac{C\lambda}{4}(2\partial\bar{\partial}\psi
-C\partial\psi\bar{\partial}\psi)\\
& \ge & -\frac{C\lambda}{4}\partial\bar{\partial}\psi =
-\frac{C\lambda}{4}g_{1/l}.
\end{eqnarray*}
Hence by (7), we find a constant $C'>0$ independent of $X_0$ such
that
\begin{equation}
\partial\bar{\partial}\lambda \ge C'\partial\bar{\partial} |z|^2,\
\ \ {\rm on\ }\Delta^{3g-3}.
\end{equation}

Let us recall the following well-known $L^2-$estimate:

\bigskip

{\bf Theorem.} (cf. [2], [9]) {\em Let $M$ be a complete K\"ahler
manifold of dimension $m$ and
let $\varphi $ be a $C^\infty $ strictly psh function on $M$. Then for any $%
\bar{\partial}-$closed $(m,1)$ form $w$ with $\int_M|w|_{\partial
\bar
\partial \varphi }^2e^{-\varphi}dV_\varphi<\infty$, there is an
m-form $u$ on $M$ such that $\bar \partial u=w$ and
\[
\left| \int_Mu\wedge \bar ue^{-\varphi }\right| \le
\int_M|w|_{\partial \bar
\partial \varphi }^2e^{-\varphi }dV_\varphi
\]
where $dV_\varphi$ denotes the volume with respect to $\partial\bar{\partial}%
\varphi$.}

\bigskip

Let $\chi:{\bf R}\rightarrow [0,1]$ be a smooth function such that
$\chi|_{(-\infty,1/2]}=1$ and $\chi|_{[1,+\infty)}=0$. Set
$w=z^\alpha \bar{\partial}\chi(|z|)\wedge dz_1\wedge\cdots\wedge
dz_{3g-3}$. Applying the above theorem for $M={\rm Teich}(R)$ with
respect to the regularization of the following psh function
$\varphi$ from above (cf. [11])
$$
\varphi=N\lambda+2(3g-3+|\alpha|)\chi(|z|)\log |z|
$$
for sufficiently large constant $N$, we obtain a form $u$ on
Teich$(R)$ satisfying $\bar{\partial}u=w$ and
$$
\left|\int_{{\rm Teich}(R)} u\wedge \bar{u}
e^{-\varphi}\right|=O(1)
$$
because of (7), (8). Then we obtain a holomorphic $3g-3$ form $s$
on Teich$(R)$ by setting $s=z^\alpha
\chi(|z|)dz_1\wedge\cdots\wedge dz_{3g-3}-u$ such that
\begin{equation}
\frac{\partial^\alpha s^\ast}{\partial z^\alpha}(0)=1,\
\frac{\partial^\beta s^\ast}{\partial z^\beta}(0)=0,\ \forall\,
|\beta|<|\alpha|,\ {\rm and \ } \|s\|_2=O(1)
\end{equation}
since $\varphi<0$ and $\varphi\sim 2(3g-3+|\alpha|)\log|z|$ near
$0$. By (4), (5), (6) and (9), the proof of Theorem 1.1 is
complete.

\bigskip

Before proving Theorem 1.2, let us recall the following

\bigskip

\textbf{Definition.} (cf. [1]) Suppose that $(M,\omega)$ is a
complete K\"ahler manifold of dimension $m$. We say that
$(M,\omega)$ has \emph{bounded geometry} if and only if for each
$x_0\in M$ there exists an embedded polydisk
$$
\iota:(\Delta^m,0)\rightarrow (M,x_0)
$$
such that the Euclidean metric and $\omega$ are equivalent on
$\Delta^m$ and for any integer $l$, there is a constant $C_l>0$
such that for any multi-indices $\alpha,\beta$ with
$|\alpha|+|\beta|\le l$ we have
$$
\left|\frac{\partial^{|\alpha|+|\beta|}}{\partial z^\alpha\partial
z^{\bar{\beta}}} g_{i\bar{j}}\right|\le C_l
$$
on $\Delta^m$ where $\omega=\sum g_{i\bar{j}}dz_i dz_j$.

\bigskip

By extreme properties of the derivatives of the Bergman metric
similar as (4), (5), it is not difficult to verify that the
Teichm\"uller space equipped with the Bergman metric has bounded
geometry. According to the Schwarz lemma of Yau (cf. Theorem 3 in
[17]), one has
$$
\frac{dV_{KE}}{dV_B}=O(1)
$$
where $dV_{KE}$ and $dV_B$ denote the volume forms of the
K\"ahler-Einstein and Bergman metrics respectively. Now we view
Teich$(R)$ as a bounded domain in ${\bf C}^{3g-3}$ equipped with
the canonical coordinate $\zeta$. If we write
$$
dV_{KE}= V_{KE}(\zeta) (\partial\bar{\partial}|\zeta|^2)^{3g-3},
$$
then
$$
V_{KE}\ge C (\delta_T|\log \delta_T|)^{-2}
$$
(cf. [8]). Note that
\begin{eqnarray*}
K_{{\rm Teich}(R)} & = & K_{{\rm Teich}(R)}^\ast dz_1\wedge \cdots
\wedge dz_{3g-3}\wedge d\bar{z}_1\wedge \cdots \wedge
d\bar{z}_{3g-3}\\
& = & K_T d\zeta_1\wedge \cdots \wedge d\zeta_{3g-3}\wedge
d\bar{\zeta}_1\wedge \cdots \wedge d\bar{\zeta}_{3g-3},
\end{eqnarray*}
which implies
$$
K_T=K_{{\rm Teich}(R)}^\ast \cdot |{\rm det}(\partial z_j/\partial
\zeta_k )|^2.
$$
Since the Bergman and Teichm\"uller metrics are equivalent, one
has
\begin{eqnarray*}
\frac{d V_B}{(\partial\bar{\partial} |\zeta|^2)^{3g-3}} & \asymp &
|{\rm det}(\partial z_j/\partial \zeta_k )|^2  \asymp  K_T.
\end{eqnarray*}
Hence $V_{KE}=O(K_T)$, verifying $(i)$ of Theorem 1.2.

\bigskip

Since ${\rm det}(\partial z_j/\partial \zeta_k )$ is nowhere
vanishing on $\Delta^{3g-3}$, we can take a single-valued branch
of $f$ of $\log {\rm det}(\partial z_j/\partial \zeta_k )$.
Applying the Schwarz-Pick lemma to the holomorphic map
$f:\Delta^{3g-3}\rightarrow \{w\in {\bf C}:|{\rm Im\,}w|<\pi\}$,
we obtain
$$
\|\partial f\|_{\partial \bar{\partial}|z|^2}(X_0)=O(1),
$$
which implies
\begin{eqnarray}
&& \|\partial \log K_T\|_{\partial \bar{\partial}\log
K_T}(X_0)\nonumber\\
& = & O (\|\partial \log K_{{\rm Teich}(R)}^\ast\|_{\partial
\bar{\partial }|z|^2}(X_0)+\|\partial
f\|_{\partial\bar{\partial}|z|^2}(X_0))= O(1).
\end{eqnarray}
The assertion $(ii)$ then follows from $(10)$ and $(i)$.

\bigskip

\textbf{Remark.} By (10), the function $r=-e^{-\tau \log
K_T}=-K_T^{-\tau}$ is a bounded strictly psh exhaustion function
on Teich$(R)$ for sufficiently small $\tau>0$. Clearly
$$
\frac1C \delta_T^{c_1}\le -r\le C \delta_T^{c_2}
$$
holds for suitable $C,c_1,c_2>0$, since trivially one has
$K_T=O(\delta_T^{-6g+6})$. Some bounded psh exhaustion functions
without estimate were given in [4], [18].

\end{document}